\documentclass[12pt]{amsart}
\def\ds{\displaystyle}
\def\CC{{\mathbb C}}

\def\NN{{\mathbb N}}
\def\OO{{\mathcal O}}
\def\ord{\hbox{ord}}
\def\phi{\varphi}

\renewcommand{\phi}{\varphi}

\title{On the product property for the Lempert function}

\author{Nikolai Nikolov and  W\l odzimierz Zwonek}

\begin{document}

\footnote{


The second-named author was supported by the KBN grant no. 5 P03A 033 21.}

\begin{abstract}
We study the problem of the product property for the Lempert
function with many poles and consider some properties of this
function mostly for plane domains.
\end{abstract}

\maketitle

\section{Introduction}

Let $A$ and $B$ be at most countable non-empty subsets of domains
$D$ and $G$ in $\Bbb C^n$ and $\Bbb C^m,$ respectively. We say
that the Lempert function $l_{D\times G}(A\times B,\cdot)$ with
pole set $A\times B$ {\it has the product property at the point
$(z,w)\in D\times G$} if $$l_{D\times G}(A\times
B,(z,w))=\max\{l_D(A,z),l_G(B,w)\}.$$ It is easy to see that this
property is true if $A$ and $B$ are singletons (cf.
\cite{Jar-Pfl}). Moreover, a necessary and sufficient condition
for the product property for the Lempert function has been given
in \cite{Dieu-Trao} (Theorem 4.1), when $B$ is a fixed singleton
and $A$ varies over all finite subsets of $D;$ namely, the product
property holds if and only if $l_G(B,w)$ is equal to the pluricomplex
Green function $g_G(B,w)$ with pole at $B.$ Unfortunately, the
proof of this result contains a gap (more precisely, there is a
gap in the proof of Lemma 2.3). A main purpose of this paper is to
prove a more general version of this lemma (Lemma 4 below), which
allows us not only to give a corrected proof of the above
mentioned result but also to refine it (Theorem 5).

Concerning the case when the pole sets are not singletons, it has
been shown in \cite{Dieu-Trao} that the product property for the
Lempert function is not true even in the case of the unit bidisc
$\Bbb D^2.$ So it is natural to study when this property holds for
$\Bbb D^2$ when each of the pole sets $A$ and $B$ has two
elements. A second purpose of the paper is to show that if, in
addition, $l_{\Bbb D}(A,0)=l_{\Bbb D}(B,0)>0,$ then the product
property for $l_{\Bbb D^2}(A\times B,(0,0))$ is true if and only
if there is a rotation sending $A$ to $B$ (Theorem 7). This result
allows us to construct easy various examples of arbitrarily large
pole sets of the unit disc for which the product property for the
Lempert function of the bidisc is not satisfied.

The paper is organized as follows. In Section 2 we give basic
facts about the Lempert function and its variations. In Section 3
we obtain explicit formulas for these functions in the plane case
and descriptions of their extremal discs (which may be considered
as an analogue of geodesic curves). These results are used in
Section 5 to construct various counterexamples to the product
property of the Lempert function. Section 4 contains proofs of
Lemma 4 and Theorem 5 mentioned above.

\section{Preliminaries}

Let $D$ be a domain in $\Bbb C^n.$ Let $z\in D$ and let $A$ be at
most countable non-empty subset of $D$ (in the paper we consider
only such sets). Denote by $\Bbb D$ the unit disc in $\Bbb C$ and
define
$$l_D(A,z):=\inf\{\prod\sb{a\in\Bbb D}|\lambda_a|\},$$ where the infimum
is taken over all subsets $(\lambda_a)_{a\in A}$ of $\Bbb D$ for
which there exists a $\phi\in\OO(\Bbb D,D)$ with $\phi(0)=z$ and
$\phi(\lambda_a)=a$ for any $a\in A$ (it is shown in
\cite{Nik-Pfl} that there are such subsets). The function
$l_D(A,\cdot)$ is called the {\it Lempert function with poles at
$A$} (cf. \cite{Dieu-Trao,Edi,Jar-Pfl1,Nik-Pfl,Wik,Wik1}). Note
that $l_D(a,\cdot):=l_D(\{a\},\cdot)$ is the classical Lempert
function. The Lempert function is monotone under inclusion of pole
sets; moreover (see \cite{Nik-Pfl}),
$$l_D(A,z)=\inf\{l_D(B,z):B\hbox{ is a finite non-empty subset of }A\},$$
and therefore,
$$l_D(A,z)=\inf\{l_D(B,z):\emptyset\neq B\subset A\}.$$

For any fixed $N\in\Bbb N^*:=\Bbb N\cup\{\infty\}$ and $a,z\in D,$
set (see \cite{Dieu-Trao})
$$l_D^N(a,z):=\inf\{\prod\sb{j=1}^N|\lambda_j|\},$$ where the
infimum is taken over all subsets $(\lambda_j)_{j=1}^N$ of $\Bbb
D$ for which there exists a $\phi\in\OO(\Bbb D,D)$ with
$\phi(0)=z$ and $\phi(\lambda_j)=a$, $j=1,\ldots,N$ (obviously,
there are such subsets). Note that we may also define another
function, denote it by $\hat l_D^N(a,z)$ in a similar way as above
but we allow some of the $\lambda_j$'s to be equal and we count
them not more than the multiplicity $\ord_{\lambda_j}\phi$ of
$\phi$ at $\lambda_j$. We shall show that both functions coincide.

{\bf Claim.} $\hat l_D^N(a,z)=l_D^N(a,z).$

{\bf Proof.} We follow the ideas of A. Edigarian (see e.g.
\cite{Edi}) who shows the result for $N=\infty$. Assume that
$N<\infty$.

It is sufficient to get that $\hat l_D^N(a,z)\geq l_D^N(a,z).$ Let
$\epsilon>0$ be arbitrary. Then we may find $\phi\in\OO(\Bbb D,D)$
such that $\phi(0)=z$ and
$\phi(\xi)-a=\prod_{j=1}^l(\xi-\lambda_j)^{k_j}\psi(\xi),$ where
$\lambda_1,\dots,\lambda_l$ are pairwise distinct numbers with
$\sum_{j=1}^l k_j=N,$ $k_j\ge 1,$ and
$$\prod_{j=1}^l|\lambda_j^{k_j}|\le\hat l_D^N(a,z)+\epsilon.$$
Consider the mapping $\phi_t(z)=\phi(t z),$ $1>t>\max_{1\le j\le
l}|\lambda_j|.$ Then $\phi_t(\Bbb D)\subset\subset D,$
$\phi_t(0)=z$ and
$$\phi_t(\xi)-a=\prod_{j=1}^l(\xi-\frac{\lambda_j}{t})^{k_j}\psi_t(z).$$
Since $\psi_t$ is a bounded, it follows that for any $s=s_t<1,$
sufficiently close to 1, the mapping $\phi_{s,t}$ defined by the
formula
$$\phi_{s,t}(\xi)=a+\prod_{j=1}^l\prod_{m=1}^{k_j}(\frac{\xi}{s^m}-\frac{\lambda_j}{t})
\psi_t(\xi)$$ belongs to the family $\OO(\Bbb D,D),$
$\phi_{s,t}(0)=z,$ and the zeroes of the double product are
pair-wise different. Thus
$$l_D^N(a,z)\leq \prod_{j=1}^l\prod_{m=1}^{k_j}s^m\frac{\lambda_j}{t}.$$
Letting $t\to 1,$ $\epsilon\to 0$, $s\to 1,$ we complete the proof.\qed

Denote by $g_D(A,\cdot)$ {\it the pluricomplex Green function with
pole at $A\subset D$}, i.e. $$g_D(A,z):=\sup\{\exp(u(z))\},$$
where the supremum is taken over all $u:D\mapsto[-\infty,0)$ such
that $u(\cdot)-\log||\cdot-a||$ is bounded from above near any
$a\in A.$ Is it known that (cf. \cite{Edi}) $$g_D(A,z)=
\inf\{\prod\sb{\lambda\in\Bbb D}\chi_A(\phi(\lambda))|\lambda|\}=
\inf\{\prod\sb{\lambda\in\Bbb D}\chi_A(\phi(\lambda))\ord_\lambda
\phi|\lambda|\} ,$$ where the infimum is taken over all
$\phi\in\OO(\Bbb D,D)$ with $\phi(0)=z.$ In particular,
$g_D(A,z)\le\tilde l_D(A,z)$ and the pluripotential Green function
$g_D(a,z):=g_D(\{a\},z)$ is equal to $\inf_{N\in\Bbb N}
l_D^N(a,z)$ (cf. \cite{Edi}). We shall see bellow that this
infimum coincides with $l_D^\infty(a,z).$

Note that the function $g_D(A,\cdot)$ is plurisubharmonic (cf.
\cite{Edi}), and the functions $l_D(A,\cdot),$ and $l_D^N$ are
upper semicontinuous.

{\bf Proposition 1.} The sequence $(l_D^N(a,z))_{N\in\Bbb N}$ is
decreasing and converges to $l_D^\infty(a,z).$

{\bf Proof.} Without loss of generality assume that $a\neq z$. To
see the inequality
$$l_D^N(a,z)\ge l_D^{N+1}(a,z),$$ let $\phi\in\OO(\Bbb D,D)$ be a
competitor for $l_D^N(a,z)$ and $\lambda_1,\dots,\lambda_N\in\Bbb
D$ be preimages of $z.$  For $\alpha\in\Bbb D,$ denote by
$$\Phi_{\alpha}(z):=\frac{\alpha-z}{1-\bar \alpha z}$$
the M\"obius transformation. Observe that if $|\alpha|<1$ is close
to $1,$ then one of the roots of the equation
$z\Phi_{\alpha}(z)=\lambda_j,$ say $\lambda_{j,1}$ is close to
$\lambda_j$ and the other one to $1$. Then it is sufficient to
take $\phi(z\Phi_{\alpha}(z))\in\OO(\Bbb D,D)$ as a competitor for
$l_D^{N+1}(a,z)$ and
$\lambda_{1,1},\lambda_{2,1},\dots,\lambda_{N,1},\lambda_{N,2}\in\Bbb
D$ as preimages of $z.$

To show that $$\lim_{N\to\infty}l_D^N(a,z)=l_D^\infty(a,z),$$
observe first that $$\liminf_{N\to\infty}l_D^N(a,z)\le
l_D^\infty(a,z)$$ by definitions. Thus, we have to prove that
$$l_D^N(a,z)\ge l_D^\infty(a,z)$$ for any $N\in\Bbb N.$ Set
$f(z)=z\exp(\frac{z-1}{z+1}).$ We claim that for every
$\lambda\in\Bbb D\setminus\{0\}$ there are infinitely many
solutions of the equation $f(z)=\lambda$ from $\Bbb D$ and that
the product of absolute values of these solutions coincides with
$|\lambda|$. Indeed, this follows from the fact that the function
$\Phi_\lambda\circ f$ has no zero radial limits and hence it is an
infinite Blaschke product (cf. \cite{Edi}). To complete the proof,
similarly as above, we consider compositions of $f$ with the
competitors for $l_D^N(a,z).$ \qed

\section{Explicit formulas and extremal discs for $l_D(A,z)$
and $l_D^N(a,z)$}

To obtain counterexamples to the product property of the Lempert
functions, we shall need explicit formulas for $l_D(A,z)$ and
$l_D^N(a,z),$ and descriptions of the extremal discs for these
functions in the plane case.

A mapping $\phi\in\OO(\Bbb D,D)$ is called an {\it
$l_D(A,z)$-extremal disc} if $\phi(0)=z$ and there exists a
nonempty subset $B$ of $A$ with $l_D(A,z)=\prod_{a\in B}|\lambda_a|,$
where $\phi(\lambda_a)=a$ for any $a\in B$ (cf.
\cite{Wik,Wik1}). A mapping $\phi\in\OO(\Bbb D,D)$ is said to be
an {\it $l_D^N(a,z)$-extremal} if $\phi(0)=z$ and
$l_D^N(a,z)=\prod\sb{j=1}\sp{M}|\lambda_j|$, where $1\le M\le N,$
$\phi(\lambda_j)=a$, $j=1,\ldots,M$, and we allow some of the
$\lambda_j$'s to be equal but they cannot be counted more than the
multiplicity $\ord_{\lambda_j}\phi$ of $\phi$ at $\lambda_j$
(compare with the definition of $\hat l_D^N$).

It is an easy observation that if $D$ is taut (i.e, if the family
$\OO(\Bbb D,D)$ is normal), $\emptyset\neq A\subset D$ is finite,
$a\in D$ and $N\in\Bbb N,$ then there are $l_D(A,z)$-extremal discs
and $l_D^N(a,z)$-extremal discs. Moreover, in this case the functions
$l_D(A,\cdot)$ and $l_D^N(a,\cdot)$ are continuous.

Recall that a plane domain $D$ is taut if and only if its boundary
contains more than one point (cf. \cite{Jar-Pfl}). On the other
hand, if the boundary of a plane domain $D$ contains at most one
point, then $l_D(A,\cdot)\equiv 0$ and $l_D^N(a,\cdot)\equiv 0.$

An application of the Schwarz-Pick lemma gives us the following
explicit formulas in the case of the unit disc: $$l_{\Bbb
D}(A,z)=\tilde l_{\Bbb D}(A,z)=\prod_{a\in
A}\bigl|\Phi_a(z)\bigr|,\ l_{\Bbb
D}^N(a,z)=\bigl|\Phi_a(z)\bigr|.$$ Moreover, if $z\not\in A,$ then
the $l_{\Bbb D}(A,z)$-extremal discs are the automorphisms
of $\Bbb D,$ sending $0$ to $z$. If $z\neq a,$ $n\in\Bbb N$, then
one may easily see that the $l_{\Bbb D}^N(a,z)$-extremal
discs are the Blaschke products of degree less than or equal to
$N,$ which map $0$ into $z.$

Now, we are going to deal with the non-simply connected plane
domains whose boundaries contain more than one point.

{\bf Proposition 2.} {\it Let $D$ be a non-simply connected plane
domain whose boundary contains more than one point, $a,z\in D$ and
$N\in\Bbb N^\ast$. Let $\pi\in\OO(\Bbb D,D)$ be a cover map with
$\pi(0)=z$. Assume that $\pi^{-1}(a)=\{\eta_1,\eta_2,\dots\}$ and
$|\eta_1|\le|\eta_2|\le\dots.$ Then
$$l_D^N(a,z)=\prod_{j=1}^N|\eta_j|.\eqno{(1)}$$

In particular, $l^N_D(a,z)<l^K_D(a,z)$ for $z\neq a$ and
$K<N\le\infty.$

If $A\subset D,$ then $$l_D(A,z)=\prod_{a\in A}\min\{|\eta|:\pi(\eta)=a\}.$$

Moreover, if $z\neq a$ and $N\in\NN$ (respectively, $l_D(A,z)>0$),
then the $l_D^N(a,z)$-extremal discs (respectively, $l_D(A,z)$-extremal discs)
are the functions of the form $\pi\circ r$, where $r$ is a rotation.}

{\bf Proof.} We shall only prove the statements for $l_D^N(a,z),$
since the proof for $l_D(A,z)$ is similar.

Without loss of generality we may assume that $a\neq z$. Let
$\varphi\in\OO(\Bbb D,D)$ be an $l_D^N(a,z)$-extremal disc. Note
that there exists $r\in\OO(\Bbb D,\Bbb D)$ with $r(0)=0$ and
$\varphi=\pi\circ r$. Choose sequences $(\lambda_{j,k})\subset\Bbb
D$, $(\nu_{j,k})\subset\Bbb N$, where $l_j\geq 1$, $j=1,\ldots,M$,
$k=1,\ldots,l_j$ and $(\lambda_j)\subset\Bbb D$, $j=1,\ldots,M$
such that all $\lambda_{j,k}$ (and all $\lambda_j$) are pairwise
different $\sum\sb{j=1}\sp{M}\sum_{k=1}^{l_j}\nu_{j,k}\leq N$,
$\nu_{j,k}\leq\ord_{\lambda_{j,k}}\phi$,
$r(\lambda_{j,k})=\lambda_j$ and
$$
\prod_{j=1}^{M}\prod_{k=1}^{l_j}|\lambda_{j,k}|^{\nu_{j,k}}=l_D^N(a,z).
$$
Certainly, $M\leq N$. Note that
$\ord_{\lambda_{j,k}}\phi=\ord_{\lambda_{j,k}}r$. Then it easily
follows from the Schwarz Lemma that
$$
\prod_{k=1}^{l_j}|\lambda_{j,k}|^{\nu_{j,k}}\geq|\lambda_j|,\
j=1,\ldots,M.
$$
Therefore,
$$
\prod_{j=1}^{M}|\lambda_j|\leq\prod_{j=1}^M
\prod_{k=1}^{l_j}|\lambda_{j,k}|^{\nu_{j,k}}\leq\prod_{j=1}^N|\mu_j|.
$$
Since $\pi(\lambda_j)=z$, $j=1,\ldots,M$, we easily get from the
way we chose $\mu_j$ that $M=N$ and, up to a permutation of the
sequence $(\lambda_j)$ we also have $|\lambda_j|=|\mu_j|$,
$j=1,\ldots,N$ and the inequalities above become equalities, which
in view of the Schwarz Lemma implies that $l_j=\nu_{j,1}=1$,
$j=1,\ldots,N$ and finally $r$ is a rotation.

 \qed

Recall now that if the boundary of a plane domain is a polar set,
then the usual Green function vanishes identically (see e.g.
\cite{Jar-Pfl}). Otherwise, we have the following description of
the $l_D^{\infty}$-extremal discs (see \cite{Nik-Zwo}), which
completes the picture in the plane case.

{\bf Proposition 3.} {\it Let $D$ be an arbitrary plane domain
whose boundary is a non-polar set, $z\in D,$ and let
$\pi\in\OO(\Bbb D, D)$ be a cover map with $\pi(0)=z.$ If $z\neq
a,$ then the $l_D^{\infty}(a,z)$-extremal discs exist and they have the
form $\pi\circ B,$ where $B\in\OO(\Bbb D,D),$ $B(0)=0$ and
$\Phi_\eta\circ B$ is a Blaschke product for any
$\eta\in\pi^{-1}(a)$.}

\section{Product property of the Lempert function}

It is known that the Green function has the product property (cf.
\cite{Edi}). In this paragraph we shall prove a result describing
when the product property of the Lempert function holds if one the
pole sets is singleton. It is a slight generalization of Theorem
2.1 in \cite{Dieu-Trao}. As we already mentioned, the main point
in the proof will be the following lemma whose proof (of a less
general version) in \cite{Dieu-Trao} (see Lemma 2.3 there) seems
to be false.

{\bf Lemma 4.} {\it Let $N\in\Bbb N^*,$ $\mu_1,\mu_2,\dots\in\Bbb
D,$ $p=\prod_{j=1}^N|\mu_j|,$ and $q\in(p,1).$ Then there exist
$f\in\OO(\Bbb D,\Bbb D)$ and $\eta_1,\eta_2,\dots\in\Bbb D$ such
that $\prod_{j=1}^N|\eta_j|=q,$ $f(0)=0$ and $f(\eta_j)=\mu_j$ for
any $j.$}

{\bf Proof.} The proof is similar to that of Proposition 1.

We may assume that $\mu_j\neq 0$ for any $j$. Otherwise we replace
the numbers $\mu_j$ by their non-zero preimages under the mapping
$z\Phi_{\alpha}(z)\in\OO(\Bbb D,\Bbb D)$ for $\alpha\in\Bbb D$
sufficiently close to $1$. We shall consider two cases.

Let first $p\neq 0.$ We shall choose the desired function of the
form $f_a(z)=z\Phi_a(z),\ a\in[0,1).$ Note that the equation
$f_a(z)=\mu_j$ has exactly two (counted with multiplicity) roots
$z_j(a),w_j(a)$, and they both belong to $\Bbb D.$ Assuming
$|z_j(a)|\le|w_j(a)|,$ we have
$|z_j(a)|\le\sqrt{|\mu_j|}\le|w_j(a)|.$ Moreover, $|z_j(a)|$ and
$|w_j(a)|$ depend continuously on $a$ (to see this, use, for
example, the formula for the solution of the equation
$f_a(z)=\mu_j$). Note also that $|z_j(0)|=|w_j(0)|=\sqrt{|\mu_j|}$
and $$\lim_{a\to 1}|z_j(a)|=|\mu_j|, \lim_{a\to 1}|w_j(a)|=1.$$

Set $$g(a)=\prod_{j=1}^{N}|z_j(a)|,\
h(a)=\prod_{j=1}^{N}|w_j(a)|,\ a\in[0,1).$$ We claim that the
functions $g$ and $h$ are continuous and if $a\to 1,$ then
$g(a)\to p$ and $h(a)\to 1.$ We also have the equality
$g(0)=h(0)=\sqrt{p}$. The only problem with these properties is
the continuity of functions $h$ and $g$ in the case $N=\infty$, so
assume that $N=\infty$. To prove the continuity, we easily see
that both functions are upper semicontinuous. On the other hand,
their lower semicontinuity follows by the inequalities
$$g(a)\ge\prod_{j=1}^M|z_j(a)|\prod_{j=M+1}^\infty|\mu_j|,\
h(a)\ge\prod_{j=1}^M|w_j(a)|\prod_{j=M+1}^\infty|\mu_j|,\ M\in\Bbb
N,$$ and the continuity of the first products.

So, if $q\le{\sqrt p},$ then there exists an $a\in[0,1)$ with
$\prod_{j=1}^N|z_j(a)|=q$; otherwise, we find an $a\in[0,1)$ with
$\prod_{j=1}^N|w_j(a)|=q$.

We shall now consider the case $p=0.$ Having in mind the case
proved above, it is sufficient to show that there exist
$f\in\OO(\Bbb D,\Bbb D)$ and points $\eta_1,\eta_2,\dots\in\Bbb D$
such that $\ds\prod_{j=1}^\infty|\eta_j|\in(0,q),$ $f(0)=0$ and
$f(\eta_j)=\mu_j,$ for any $j.$ Fix $k$ with
$\ds\prod_{j=1}^k|\eta_j|<q^2,$ choose $\varepsilon$ such that
$$\max_{|z|\le\sqrt{|\mu_j|}}|e^{\varepsilon\frac{z-1}{z+1}}-1|
<1-\sqrt{|\mu_j|}$$ for $1\le j\le k,$ and set
$f=z\exp(\varepsilon\frac{z-1}{z+1}).$ It follows by the Rouch\'e
theorem that the functions $z-\mu_j$ and $f-\mu_j$ have the same
numbers of zeroes inside the disc $\{z\in\Bbb
C:|z|<\sqrt{|\mu_j|}\}.$ Hence for any $j\le k$ there is a unique
$\eta_j$ from this disc such that $f(\eta_j)=\mu_j.$ On the other
hand, similarly as in the proof of Proposition 1, the function
$\Phi_{\mu_j}\circ f$ is an infinite Blaschke product. Therefore,
for any $j>k,$ we may choose $\eta_j$ with
$\ds|\eta_j|>2^{-2^{-j}}$ and $f(\eta_j)=\mu_j.$ Thus
$\ds\prod_{j=1}^\infty|\eta_j|$ is a non-zero product, smaller
than $q.$\qed

Now we are going to the main result in this section.

{\bf Theorem 5.} {\it Let $D$ and $G$ be domains in $\CC^n$ and
$\CC^m,$ respectively, and let $z\in D,$ $w,b\in G.$ Then for any
nonempty at most countable $A\subset D$ the following inequalities
hold:
$$\max\{l_D(A,z),l_G^{\# A}(b,w)\}\le l_{D\times
G}(A\times\{b\},(z,w))\le\\ \max\{l_D(A,z),l_G(b,w)\}.$$ Moreover,
for given $N\in\NN^{\ast}$ the equality $$l_{D\times
G}(A\times\{b\},(z,w))=\max\{l_D(A,z),l_G(b,w)\}$$ holds for any
$A\subset D$ with $N$ elements if and only if
$l_G(b,w)=l_G^N(b,w).$}

{\bf Proof.} The proof is similar to that of Theorem 2.1 in
\cite{Dieu-Trao}. The left hand-side inequality follows by the
definitions. To prove the other one, let $\alpha<1$ be such that
$\alpha>\max\{l_D(A,z),l_G(b,w)\}.$ If $A=\{a_j\}_{j=1}^N,$ then
there exist $\varphi\in\OO(\Bbb D,D),$ $\lambda_j\in\Bbb D,$
$\psi\in\OO(\Bbb D,G)$ and $\zeta\in\Bbb D$ such that
$\varphi(0)=z,$ $\varphi(\lambda_j)=a_j,$ $\psi(0)=w,$
$\psi(\zeta)=b,$ and $\max\{\prod_{j=1}^N|\lambda_j|,\
|\zeta|\}<\alpha$.  By Lemma 4 we may find $f\in\OO(\Bbb D,\Bbb
D)$ and $\eta_1,\eta_2,\dots\in\Bbb D$ such that
$\prod_{j=1}^N|\eta_j|=\alpha,$ $f(0)=0$ and
$f(\eta_j)=\lambda_j.$ Set
$$B=\prod_{j=1}^N\frac{\bar\eta_j}{|\eta_j|}\Phi_{\eta_j},\
\xi=(\varphi\circ f,\psi(\frac{\zeta}{\alpha}\Phi_{\alpha}\circ
B).$$ Then $\xi\in\OO(\Bbb D,D\times G),$ $\xi(0)=(z,w)$ and
$\xi(\eta_j)=(a_j,b),$ which implies that $$l_{D\times
G}(A\times\{b\},(z,w))\le\alpha.$$ Hence $$l_{D\times
G}(A\times\{b\},(z,w))\le\max\{l_D(A,z),l_G(b,w)\}$$

It remains to show that if $$l_{D\times
G}(A\times\{b\},(z,w))=\max\{l_D(A,z),l_G(b,w)\}$$ for any
$A\subset D$ with $N$ elements, then $l_G(b,w)\le l_G^N(b,w),$ because
the opposite inequality always holds. We know that for any
$\varepsilon>0$ there exist $\varphi\in\OO(\Bbb D,G)$ and pair-wise
distinct points $\eta_1,\eta_2,\dots\in \Bbb D$ such that
$\varphi(0)=w,$ $\varphi(\eta_j)=b,\ j=1,2,\dots$ and
$$\prod_{j=1}^\infty|\eta_j|<g_G(b,w)+\varepsilon.$$ Note that we
may choose $\psi\in\OO(\Bbb D,D)$ with $\psi(0)=w$ and
$\psi(\eta_j)\neq\psi(\eta_k)$, if $j\neq k$. Set
$A=\{\psi(\eta_j)\}_{j=1}^N.$ Since $(\psi,\phi)\in\OO(\Bbb D,
D\times G)$ is a competitor for $l_{D\times
G}(A\times\{b\},(z,w)),$ we conclude that $$l_G(b,w)\le l_{D\times
G}(A\times\{b\},(z,w))\le l_G^N(b,w)+\varepsilon.\qed$$

{\bf Corollary 6.} {(see \cite{Dieu-Trao})} {\it Let $D$ and $G$ be
domains in $\CC^n$ and $\CC^m,$ respectively, and let $z\in D,$
$w,b\in G.$ Then the equality $$l_{D\times
G}(A\times\{b\},(z,w))=\max\{l_D(A,z),l_G(b,w)\}$$ holds for any
nonempty at most countable $A\subset D$ if and only if
$l_G(b,w)=g_G(b,w).$}

Note that by the Lempert theorem (cf. \cite{Jar-Pfl}) the last
equality holds for any convex domains. It is also true for the
symmetrized bidisc which is not biholomorphic to a convex domain
(see \cite{Cos}, see also \cite{Jar-Pfl2}).

\section{Counterexamples to the product property of the Lempert
function}

Let $G$ be a plane domain and let $D$ be a domain in $\Bbb C^n.$
Theorem 5 and the explicit formula for $l_G^N$ (Proposition 2)
show that the product property for $l_{D\times
G}(A\times\{b\},(z,w)),$ $b\neq w,$ holds if and only if either
$G$ is simply connected or its complement is a singleton.

In this paragraph we shall see that the product property for the
Lempert function of the bidisc is a seldom phenomenon if each of
the pole sets has more than one element.

We also show that the left-hand side inequality in Theorem 5 is
not a good candidate for a modified product property; namely, this
inequality is strict, in general, for non-simply connected domains
whose boundaries contain more than one point.

Since the Green function has the product property, it does not
exceed the Lempert function and both functions coincide on the
unit disc, it follows that $$\max\{l_{\Bbb D}(A,z),l_{\Bbb
D}(B,w)\}\le l_{\Bbb D^2}(A\times B,(z,w)).\eqno{(2)}$$ On
the other hand, we have

{\bf Theorem 7.} {\it Let $A$ and $B$ be two-point subsets of
$\Bbb D,$ such that $0\not\in A,B$ and $l_{\Bbb D}(A,0)=l_{\Bbb
D}(B,0).$ Then $$l_{\Bbb D}(A,0)=l_{\Bbb D^2}(A\times B,(0,0))$$
if and only if there is a rotation sending $A$ to $B.$

In addition, if $B=e^{i\theta}A,$ $\theta\in\Bbb R,$ then the
$l_{\Bbb D^2}(A\times B,(0,0))$-extremal discs are of the form
$(r,e^{i\theta}r),$ where $r$ is a rotation.}

{\bf Proof.} Let $A=\{a_1,a_2\},$ $B=\{b_1,b_2\},$ and let $\psi=(\psi_1,\psi_2)$
be an $l_{\Bbb D^2}(A\times B,(0,0))$-extremal disc. Then
there are a set $J\subset \{1,2\}\times\{1,2\}$ and numbers
$z_{k,l}\in\Bbb D,$ $(k,l)\in J,$ such that
$$\psi(z_{k,l})=(a_k,b_l)\hbox{ and } \prod_{(k,l)\in J}|z_{k,l}|=
l_{\Bbb D^2}(A\times B,(0,0)).$$

Suppose now that $l_{\Bbb D}(A,0)=l_{\Bbb D^2}(A\times B,(0,0)).$ Since
$$l_{\Bbb D}(A,0)=l_{\Bbb D}(a_1,0)\cdot l_{\Bbb D}(a_2,0),$$ it follows
that $\psi_1$ is an $l_{\Bbb D}^2(a_j,0)$-extremal disc, $j=1,2$
Analogously, $\psi_2$ is an $l_{\Bbb D}^2(b_j,0)$-extremal disc,
$j=1,2$. In particular, $\#J\ge 2,$ and $\psi_1,\psi_2$ are rotations or
Blaschke products of degree two.

If $\#J=3,$ then $\psi_1$ or $\psi_2$ must be simultaneously a rotation and a
Blaschke product of degree two, which is a contradiction.

Let $\#J=4.$ Then we may assume that
$$\psi_1(z)=z\Phi_\alpha(z),\ \psi_2(z)=e^{it}z\Phi_\beta$$ for
some $\alpha,\beta\in\Bbb D,t\in\Bbb R$. Therefore,
$$z_{1,1}\Phi_\alpha(z_{1,1})=z_{1,2}\Phi_\alpha(z_{1,2}),\
z_{2,1}\Phi_\alpha(z_{2,1})= z_{2,2}\Phi_\alpha(z_{2,2}),$$
$$z_{1,1}\Phi_\beta(z_{1,1})= z_{2,1}\Phi_\beta(z_{2,1}),\
z_{1,2}\Phi_\beta(z_{1,2})= z_{2,2}\Phi_\beta(z_{2,2}).$$ It
follows that that
$$z_{1,1}=\Phi_\alpha(z_{1,2})=\Phi_\beta(z_{2,1}),\
z_{1,2}=\Phi_\beta(z_{2,2}),\ z_{2,1}=\Phi_\alpha(z_{2,2}).$$
Hence $z_{1,1}=\Phi_\alpha\circ\Phi_\beta(z_{2,2})=
\Phi_\beta\circ\Phi_\alpha(z_{2,2}).$ Then a straightforward
calculation leads to the equality
$$(2-\alpha\bar\beta-\bar\alpha\beta)(z_{2,2}^2(\bar\alpha-\bar\beta)
+z_{2,2}(\alpha\bar\beta-\bar\alpha\beta)+\beta-\alpha)=0.$$ It is
easy to see that if $\alpha\neq\beta,$ then both roots of the
equation $$z^2(\bar\alpha-\bar\beta)
+z(\alpha\bar\beta-\bar\alpha\beta)=\alpha-\beta$$ belong to
$\partial\Bbb  D.$ Thus, $\alpha=\beta,$ $z_{1,2}=z_{2,1},$
$z_{1,1}=z_{2,2}$, a contradiction.\qed

It remains to consider the case $\#J=2.$
Then either $(1,2),(2,1)\not\in J,$ or $(1,1),(2,2)\not\in J.$
It follows that $\psi_1$ and $\psi_2$ must be rotations, say
$\psi_1(z)=e^{i\theta_1}z,$ $\psi_2(z)=e^{i\theta_2}z,$ $\theta_1,
\theta_2\in\Bbb R,$ and hence $B=e^{i\theta}A,$ where
$\theta=\theta_1-\theta_2.$

Conversely, it is clear that if $B=e^{i\theta}A$ and $r$ is a
rotation, then the mapping $(r,e^{i\theta}r)\in\OO(\Bbb D,\Bbb
D^2)$ is a competitor for $l_{\Bbb D^2}(A\times B,(0,0)).$
This implies that $$l_{\Bbb D}(A,0)\geq\ l_{\Bbb D^2}(A\times
B,(0,0)).$$ Now the inequality $(2)$ completes the proof.\qed

A consequence of Theorem 7 is the following

{\bf Corollary 8.} {\it Let $A,B$ be two-point subsets of $\Bbb D$ and
$z\in\Bbb D\setminus A.$ Then there exist uncountable many $w\in\Bbb D$
for which $$l_{\Bbb D}(A,z)=l_{\Bbb D}(B,w)<l_{\Bbb D^2}(A\times
B,(z,w)).$$}

{\bf Proof.} It suffices to note that there exist uncountable many
$w$'s with $l_{\Bbb D}(A,z)=l_{\Bbb D}(B,w),$ but at most two
$w$'s for which there is an automorphism of $\Bbb D,$ sending $z$
to $w$ and $A$ to $B.$\qed

We do not know whether Theorem 7 still holds for sets with equal
numbers of elements, greater than 1. However, this theorem and the next
proposition provide for given $(z,w)\in\Bbb D^2$ a large class of
counterexamples to the product property of
$l_{\Bbb D^2}(A\times B,(z,w))$ for pole sets $A$ and $B$ with
arbitrary numbers of elements, greater than 1.

{\bf Proposition 9.} {\it Let $D$ and $G$ be domains in $\CC^n$
and $\CC^m,$ respectively. Let $z\in D,\ w\in G,\ A\subset D,\
B\subset G$ and $q\in(0,1)$ be such that
$$\max\{l_D(A,z),l_G(B,w)\}=q\l_{D\times G}(A\times B,(z,w))>0.$$
Then $$\max\{l_D(A\cup A_1,z),l_G(B\cup B_1,w)\}<
l_{D\times G}((A\cup A_1)\times(B\cup B_1),(z,w))$$ for any
$A_1\subset D,\ B_1\subset G$ with $A\cap A_1=B\cap B_1=\emptyset$
and $$g_D(A_1,z)g_G(B_1,w)>q.\eqno{(3)}$$}

{\bf Proof.} It is easy to see that
$$l_{D\times G}((A\cup A_1)\times(B\cup B_1),(z,w))\ge$$
$$l_{D\times G}(A\times B,(z,w))l_{D\times G}(A\times B_1,(z,w))
l_{D\times G}(A_1\times (B\cup B_1),(z,w))$$
$$\ge l_{D\times G}(A\times B,(z,w))g_{D\times G}(A\times
B_1,(z,w))g_{D\times G}(A_1\times (B\cup B_1),(z,w))$$
$$\ge l_{D\times G}(A\times B,(z,w))g_G(B_1,w)g_D(A_1,z)$$
$$>\max\{l_D(A,z),l_G(B,w)\}\ge\max\{l_D(A\cup A_1,z),l_G(B\cup B_1,w)\}.
\qed$$

{\bf Remark.} Recall that if the boundary of a planar domain $D$
is a non-polar set, then there exists a polar set
$F\subset\partial D$ such that $\lim_{a\to a_0}g_D(a,z)=1$ for any
$a_0\in(\partial D)\setminus F$ and any $z\in D.$ Since
$$g_D(A,z)=\prod_{a\in A}g_D(a,z),$$ it follows that for a given
$q\in(0,1)$ and $N\in\Bbb N^\ast$ there is a set $A$ with $N$
elements and with $\hbox{dist}(A,a_0)<1-q$, and $g_D(A,z)>q.$ So,
we may provide the inequality $(3)$ for any planar domains whose
boundaries are non-polar.

Now we shall prove two results showing that the left-hand side
inequality in Theorem 5 is also strict for general plane domains.

{\bf Proposition 10.} {\it Let $D$ and $G$ be plane domains whose
boundaries contain more than one point, $w,b\in G,$ $w\neq b$,
$z\in D$. Assume that $G$ is non-simply connected. Then there
exists a countable set $A=\{a_1,a_2,\dots\}$ of points in $D$ such
that if $A_N=\{a_1,a_2,\dots,a_N\},$ $N\in\NN\setminus\{1\},$ then
$$l_D(A_N,z)=l_G^N(b,w)<l_{D\times G}(A_N\times\{b\},(z,w)).$$

Moreover, if the boundary of $G$ is a non-polar set, then
$$l_D(A,z)=g_G(b,w)<l_{D\times G}(A\times\{b\},(z,w)).$$}

{\bf Proof.} Since $l_D(\cdot,z)$ is a continuous function,
$l_D(z,z)=0$ and $\lim_{a\to\partial D}l_D(a,z)=1$ (which follows
by the explicit formula for $l_D(a,z)$), we may find $a_1$ with
$l_D(a_1,z)=l_G(b,w)>0$. Using similar argument, we obtain a
sequence of points $a_1,a_2,\dots\in D$ such that
$$l_D(A_N,z)=l_G^N(b,w)>0.$$ Moreover, each of these points can be
chosen in uncountable many ways. Thus, if $\pi\in\OO(\Bbb D,G)$
and $\tau\in\OO(\Bbb D,D)$ are cover maps with $\pi(0)=w$ and
$\tau(0)=z,$ then we may assume that
$$\frac{\xi_1}{\xi_2}\neq\frac{\eta}{\zeta}\hbox{ for any }
\xi_1\in\tau^{-1}(a_1),\ \xi_2\in\tau^{-1}(a_2),\
\eta,\zeta\in\pi^{-1}(b).\eqno{(4)}$$

Suppose now that for some $N\in\NN\setminus\{1\}$ we have
$$l_D(A_N,z)=l_{D\times G}(A_N\times \{b\},(z,w)).$$ Since
$D\times G$ is a taut domain, there exists an $\tilde l_{D\times
G}(A_N\times\{b\},(z,w))$-extremal disc $(\varphi,\psi).$ Then
$\varphi$ and $\psi$ must be $\tilde l_D(A_N,z)$-extremal disc and
$l_D^N(b,w)$-extremal disc, respectively. By Proposition 2, we may
assume that $\psi=\pi$ and $\varphi=\tau\circ e^{i\theta}$ for
some real $\theta.$ In particular, there are
$\eta_1\in\pi^{-1}(b)\cap e^{-i\theta}\tau^{-1}(a_1)$ and
$\eta_2\in\pi^{-1}(b)\cap e^{-i\theta}\tau^{-1}(a_2).$ A
contradiction with $(4).$

We are going to the second part of the proposition. First, we
shall show that there exists an
$l_{D\times G}(A\times\{b\},(z,w))$-extremal disc.
Let $\xi_N,$ $N\in\Bbb N,$ be an
$l_{D\times G}(A_N\times\{b\},(z,w))$-extremal disc.
Then there exist sets $J\subset\{1,2,\dots,N\}$ and
$(\lambda_{j,N})_{j\in J}\subset\Bbb D$ such that $\xi_N(0)=(z,w)$,
$\xi_N(\lambda_{j,N})=(a_j,b)$ for any $j\in J,$ and
$$l_{D\times G}(A_N\times\{b\},(z,w))=\prod_{j\in J}|\lambda_{j,N}|.$$
Put $\lambda_{j,N}=1$ for $j\not\in J.$
Passing to subsequences and applying the standard diagonal
process, we may assume that $\xi_N\to\xi\in\OO(\Bbb D,D\times G)$
uniformly on compact subsets of $\Bbb D,$
$\lim_{N\to\infty}\lambda_{j,N}=\lambda_j\in\overline{\Bbb D}$ for
any $j$ and $\xi(0)=(z,w)$, $\xi(\lambda_j)=(a_j,b)$. To prove
that $\xi$ is an $l_{D\times G}(A\times\{b\},(z,w))$-extremal disc,
suppose the contrary. It follows that
$$\prod_{j=1}^\infty|\lambda_j|\ge ql_{D\times G}(A\times\{b\},(z,w)),$$
where $q>1.$ Then for any $k$ there is $n_k$ such that
$$\prod_{j=1}^k|\lambda_{j,N}|\ge ql_{D\times G}(A\times\{b\},(z,w))$$
if $N\ge n_k.$ Since
$$\prod_{j=k+1}^\infty|\lambda_{j,N}|\ge l_D(A\setminus A_k,z),$$
we obtain that
$$l_{D\times G}(A_N\times\{b\},(z,w))\ge
ql_D(A\setminus A_k,z)l_{D\times G}(A\times\{b\},(z,w)).\eqno{(5)}$$
Note that the assumption that
the boundary of $G$ is a non-polar set implies
$$l_D(A,z)=g_G(b,w)>0.$$ Then (among others we use the fact that the sequence
$(a_k)$ has no accumulation point in $D$)
$$\lim_{N\to\infty}l_{D\times G}(A_N\times\{b\},(z,w))=
l_{D\times G}(A\times\{b\},(z,w))>0,$$
$$\lim_{k\to\infty} l_D(A\setminus A_k,z)=1.$$
It follows by $(5)$ that $1\ge q$, a contradiction.

Thus, $\xi=(\varphi,\psi)$ is an
$l_{D\times G}(A\times\{b\},(z,w))$-extremal disc. Suppose now that
$$l_D(A,z)=g_G(b,w)=l_{D\times G}(A\times\{b\},(z,w)).$$
Then $\varphi$ and $\psi$ must be a $l_D(A,z)$-extremal
disc and an $l_G^\infty(b,w)$-extremal disc, respectively. By
Propositions 2 and 3, we may assume that $\varphi=\tau$ and
$\psi=\pi\circ B,$ where $B\in\OO(\Bbb D,\Bbb D),$ $B(0)=0$ and
$\Phi_\eta\circ B$ is a Blaschke product for any
$\eta\in\pi^{-1}(b)$. Moreover, $J=\Bbb N,$
$\lambda_j\in\tau^{-1}(a_j),$ $|\lambda_j|=l_D(a_j,z)$ and
$B(\lambda_j)\in\pi^{-1}(b)$ for any $j.$ On the other hand,
$$|\lambda_1|\ge|B(\lambda_1)|\ge l_G(b,w)=l_D(a_1,z)=|\lambda_1|,$$ which
shows that $B$ is a rotation. Then, as above, we get a
contradiction with $(4),$ which completes the proof.\qed

This proof allows us to obtain also the following

{\bf Proposition 11.} {\it Let $D$ and $G$ be plane domains whose
boundaries contain more than one point, $w,b\in G,$ $w\neq b$,
$z\in D$. Assume that $G$ is non-simply connected. Then there
exists a countable subset $A$ of $D$ such that
$$\max\{l_D(A,z),g_G(b,w)\}<l_{D\times G}(A\times\{b\},(z,w)).$$}

{\bf Proof.} Choose the points $a_1$ and $a_2$ in $D$ as in the
proof of Proposition 10 and set $A_2=\{a_1,a_2\}.$ Let $q\in(0,1)$
be such that
$$l_D(A_2,z)=l_G^2(b,w)=q l_{D\times G}(A_2\times\{b\},(z,w)).$$
Now it is enough to note that for any countable subset $B\subset D$
with $l_D(B,z)>q$ we have
$$l_{D\times G}((A_2\cup B)\times\{b\},(z,w))\ge
l_{D\times G}(A_2\times\{b\},(z,w))l_D(B,z)>$$
$$l_D(A_2,z)=l_G^2(b,w)>\max\{l_D(A_2\cup B,z),g_G(b,w)\}.\qed$$

\noindent{\bf Acknowledgments.} A part of this paper was prepared
during the stay of the first-named author at the Jagiellonian
University, Krak\'ow, Poland (October, 2003). The authors would
like to thank Marek Jarnicki and Peter Pflug for drawing their
attention to the gap in the proof of Lemma 2.3 in
\cite{Dieu-Trao}.

\bigskip

\noindent
\begin{tabular}[t]{l}
{\sc Nikolai Nikolov} \\
Institute of Mathematics and Informatics,\\
Bulgarian Academy of Sciences,\\
Acad.\ G. Bonchev Str., Block 8,\\
1113 Sofia, {\sc Bulgaria}\\[6pt]
{\itshape E-mail:} nik@math.bas.bg
\end{tabular}
\hfill
\begin{tabular}[t]{l}
{\sc W\l odzimierz Zwonek}\\
Instytut Matematyki,\\
Uniwersytet Jagiello\'nski,\\
Reymonta 4, \\
30-059 Krak\'ow, {\sc Poland}\\[6pt]
{\itshape E-mail:} zwonek@im.uj.edu.pl
\end{tabular}


\begin{thebibliography}{99}

\bibitem{Cos}{\sc C. Costara,} {\itshape The symmetrized bidisc as a
counterexample to the converse of Lempert's theorem,}
{Bull. London Math. Soc.,} to appear.

\bibitem{Dieu-Trao}{\sc N. Q. Dieu, N. V. Trao,} {\itshape Product property
of certain extremal functions,} {Complex Variables,} {\bfseries 48(8)} (2003),
681--694.

\bibitem{Edi}{\sc A. Edigarian,} {\itshape Analytic discs method in
complex analysis,} {Diss.\ Math.,} {\bfseries 402} (2002), 1-56.

\bibitem{Jar-Pfl}{\sc M. Jarnicki, P. Pflug,} {\itshape Invariant
distances and metrics in complex analysis,} {De Gruyter} (1993).

\bibitem{Jar-Pfl1}{\sc M. Jarnicki, P. Pflug,} {\itshape Invariant
distances and metrics in complex analysis - revisited,} {Diss.\ Math.,} to
appear.

\bibitem{Nik-Pfl}{\sc N. Nikolov, P. Pflug,} {\itshape The multipole Lempert
function is monotone under inclusion of sets,} Preprint (2004).

\bibitem{Nik-Zwo}{\sc N. Nikolov, W. Zwonek,} {\itshape Some remarks on the
Green function and the Azukawa pseudometric,} {Monats. Math.,}
{\bfseries 142(4)} (2004), 341--350.

\bibitem{Wik}{\sc F. Wikstr\"om,} {\itshape Non-linearity of the pluricomplex
Green function,} {Proc. Amer. Math. Soc.,} {\bfseries 129(4)} (2001), 1051--1056.

\bibitem{Wik1}{\sc F. Wikstr\"om,} {\itshape Qualitative properties of
biholomorphically functions with multiple poles,} Preprint (2004).
\end{thebibliography}
\end{document}